\documentclass[a4paper,10pt,twoside]{article}
\usepackage{pst-fill,pst-grad}
\usepackage{textcomp}
\usepackage[english]{babel}
\usepackage[utf8x]{inputenc}
\usepackage{graphicx}
\usepackage{amsmath}
\usepackage{float}
\usepackage{fancyhdr}
\usepackage[matrix,arrow,curve]{xy}
\usepackage{pstricks} 
\usepackage{amsmath,amsfonts,verbatim,afterpage,theorem,euscript,mathrsfs,amssymb}
\usepackage{amsfonts}
\usepackage{amssymb}
\usepackage{array}
\usepackage{dsfont}
\usepackage{hyperref}

\newtheorem{Proposition}{Proposition}[section]
\newtheorem{Lemme}{Lemma}[section]
\newtheorem{Theoreme}{Theorem}
\newtheorem{Corollaire}{Corollary}[section]
\newtheorem{Remarque}{Remark}[section]
\title{\bf Fractional Laplacians and Nilpotent Lie Groups}
\author{Diego Chamorro\footnote{Universit\'e d'Evry, LaMME (diego.chamorro@univ-evry.fr)} and Oscar Jarr\'in\footnote{Universit\'e d'Evry, LaMME}}

\begin{document}
\maketitle
\begin{scriptsize}
\abstract{The aim of this short article is to generalize, with a slighthly different point of view, some new results concerning the fractional powers of the Laplace operator to the setting of Nilpotent Lie Groups and to study its relationship with the solutions of a partial differential equation in the spirit of the articles of Caffarelli \& Silvestre \cite{Caffarelli} and Stinga \& Torrea \cite{Stinga1}.}\\[3mm]
\end{scriptsize}

\section{Introduction}
In this article we are interested in a generalization, from the euclidean setting $\mathbb{R}^n$ to Nilpotent Lie groups, of the relationship between the solutions $u(t,x)$ of the following partial differential equation (also called \textit{extension problem}):
\begin{equation}\label{EquationInicial1}
\begin{cases}
\partial_t^2 u(t,x)+\frac{1-2s}{t}\partial_t u(t,x)+\Delta_x u(t,x)= 0&\qquad 0<s<1 \mbox{ and } t>0\\[5mm]
u(0,x)=\varphi(x),& \qquad \varphi \in \mathcal{S}(\mathbb{R}^n),
\end{cases}
\end{equation}
and the fractional powers of the Laplacian of the initial data $\varphi$:
\begin{equation}\label{FractionalLaplacian1}
(-\Delta)^{s}\varphi\qquad \mbox{with } 0<s<1.
\end{equation}
In $\mathbb{R}^n$ endowed with its natural structure, the operator $(-\Delta)^{s}$ can be defined in the Fourier level by the formula  $\widehat{(-\Delta)^{s}\varphi}(\xi)=C |\xi|^{2s}\widehat{ \varphi}(\xi)$. We will see later on how to define this operator in the framework of nilpotent Lie groups.\\

The relationship between (\ref{EquationInicial1}) and (\ref{FractionalLaplacian1}) is given by the following expression:
\begin{equation}\label{RelationShip}
\underset{t\to 0^+}{ \lim} t^{1-2s}\partial_tu(t,x)=-C(s)(-\Delta)^s\varphi(x).
\end{equation}
This formula and many of its applications in the theory of partial differential equations was first studied by Caffarelli \& Silvestre \cite{Caffarelli} in 2007 and since then this work has been generalized to many different frameworks.\\ 

This was done for example for linear second order partial differential operators in the article of Stinga \& Torrea \cite{Stinga1}, see also Ferrari \& Franchi \cite{Ferrari} for a generalization to the setting of Carnot groups and more recently for a treatment in non-compact manifolds see Banica \& al. \cite{Banica}. See also Frank \& al. \cite{Frank}, Gal\'e \& al. \cite{Gale} and Cabr\'e \& Sire \cite{Cabre} for other interesting applications and generalizations of this work.\\

Our aim in this article is to generalize the relationship between (\ref{EquationInicial1}) and (\ref{FractionalLaplacian1}) to the setting of nilpotent Lie groups but we will not study here all the possible applications of the formula (\ref{RelationShip}) given in the aforementioned articles as this task will take a considerable amount of work and deserves further investigations.\\

In the recent article of Ferrari \& Franchi \cite{Ferrari} this relationship is studied in great details in the setting of the Carnot groups and we want to give here a slightly different point of view which is essentially based on the fact that there is not a unique way to define a Laplace operator in the framework of Lie groups. \\

Let us introduce a typical example of this situation. The euclidean setting of $\mathbb{R}^n$ has a generalization in the direction of Lie groups and perhaps one of the most well-known generalization in this sense is given by the Heisenberg group $\mathbb{H}$ which can be defined in the following manner: consider $x=(x_{1},x_{2},x_{3})$ an element of $\mathbb{R}^{3}$ and the non commutative group law
\begin{equation*}
x\cdot y=(x_{1},x_{2},x_{3})\cdot(y_{1},y_{2},y_{3})=\big(x_{1}+y_{1},x_{2}+y_{2},x_{3}+y_{3}+\frac{1}{2}(x_{1}y_{2}-y_{1}x_{2})\big).
\end{equation*}
Define next a dilation structure $\delta_\alpha$ for $\alpha>0$ by
\begin{eqnarray*}
\delta_{\alpha}: \mathbb{R}^{3}&\longrightarrow &\mathbb{R}^{3}\\
x=(x_{1},x_{2},x_{3}) &\longmapsto & \delta_{\alpha}[x]=(\alpha x_{1}, \alpha x_{2},\alpha^{2} x_{3})\nonumber
\end{eqnarray*}
The triplet $\mathbb{H}=(\mathbb{R}^3, \cdot, \delta)$ corresponds to the \textit{Heisenberg group}\footnote{For more details concerning the Heisenberg group, see \cite{Stein2}, Chapters XII and XIII.}. We can define a norm and a distance in this group by the formulas $\|x\|=\left[\left(x_{1}^{2}+x_{2}^{2}\right)^{2}+16 x_{3}^{2}\right]^{\frac{1}{4}}$ and $d(x,y)=\|y^{-1}\cdot x\|$, for all $x,y\in \mathbb{H}$. From the point of view of measure theory, the Haar measure is given by the Lebesgue measure and the volume of open balls $V(r)=V(B(x,r))$ satisfy the homogeneous identity
$$V(r)=r^4V(1).$$
Remark that the homogeneous dimension with respect to the dilation structure of the Heisenberg group is $N=4$ while its topological dimension is $n=3$.\\

Associated to this group we have a Lie algebra $\mathfrak{h}$ given by the left-invariant vector fields
$$X_{1}=\frac{\partial}{\partial x_{1}}-\frac{1}{2}x_{2}\frac{\partial}{\partial x_{3}}, \quad X_{2}=\frac{\partial}{\partial x_{2}}+\frac{1}{2} x_{1} \frac{\partial}{\partial x_{3}}\qquad \mbox{and} \qquad T=\frac{\partial}{\partial x_{3}}$$
and we have the identities
$$[X_{1},X_{2}]=X_{1}X_{2}-X_{2}X_{1}= T,\qquad  [X_{i},T]=[T,X_{i}]=0 \quad  \mbox{where } i=1,2.$$
We remark from these commutation properties that the vector fields $X_1$ and $X_2$ \textit{span} the Lie algebra $\mathfrak{h}$ of the Heisenberg group. We will say then that $X_1$ and $X_2$ form the \textit{first layer} of the stratification of the Lie algebra $\mathfrak{h}$ while $T$ lies in the \textit{second layer} of the stratification of $\mathfrak{h}$ and since all the other commutation brackets are null, we say that the Heisenberg group is a \textit{nilpotent} Lie group with a \textit{stratification} of order $2$.\\

We will say that a function on $ \mathbb{H}\setminus\{0\}$ is \textit{homogeneous} of degree $\lambda \in \mathbb{R}$ if $f(\delta_{\alpha}[x]) =  \alpha^\lambda f (x)$ for all $\alpha > 0$ and, in the same manner, we will say that a differential operator $D$ is homogeneous of degree $\lambda$ if $D\big(f(\delta_{\alpha}[x])\big) = \alpha^\lambda \big(D f \big) (\delta_{\alpha}[x])$. With these definitions we can see that the vector fields of the Lie algebra $\mathfrak{h}$ satisfy the following homogeneous properties
$$X_1\big(f (\delta_{\alpha}[x])\big) = \alpha (X_1 f )(\delta_{\alpha}[x]), \qquad X_2\big(f( \delta_{\alpha}[x])\big) = \alpha (X_2 f )(\delta_{\alpha}[x]),$$
$$\mbox{and} \qquad  T\big(f (\delta_{\alpha}[x])\big) = \alpha^2 (T f )(\delta_{\alpha}[x]).$$
We observe from this that the vector fields $X_1$ and $X_2$ are homogeneous of degree $1$ but $T$ is homogeneous of degree $2$. \\

Now we will see how these rather simple facts and properties have important consequences. Indeed, if we want to build from these vector fields an equivalent of the Laplace operator in $\mathbb{R}^n$ (which is, among other properties, a differential operator of degree 2) we have several choices:
\begin{eqnarray}
\mathcal{J}_1&=&-(X_1^2+X_2^2) \nonumber\\[3mm]
\mathcal{J}_2&=&-(X_1^2+X_2^2+T)  \label{Several_Laplacians}\\[3mm]
\mathcal{J}_3&=&-(X_1^2+X_2^2+T^2)\nonumber
\end{eqnarray}
The operator $\mathcal{J}_1$ is called the \textit{Kohn-Laplace} operator and is built solely from the vector fields of the first layer of the stratification, of course it is an operator of degree $2$ with respect to the dilation $\delta$. The second operator $\mathcal{J}_2$ is now constructed using \textit{all} the vector fields of the Lie algebra $\mathfrak{h}$, but it is also an operator of degree $2$ with respect to the dilation $\delta$. The third operator $\mathcal{J}_3$, sometimes called the \textit{full Laplacian}, is given by all the vector fields of the Lie algebra $\mathfrak{h}$, but now we can see that $\mathcal{J}_3$ is not an homogeneous differential operator with respect to the dilation structure.\\

Of course, depending on the choice of the Laplace operator, we will obtain different properties: just think for example in the heat kernel associated to it or in the underlying geometry of the group, see \cite{Furioli1} or \cite{Varopoulos} for more details.\\

The generalization of the relationship between (\ref{EquationInicial1}) and (\ref{FractionalLaplacian1}) in the setting of  stratified Lie groups was carried out by Ferrari \& Franchi in \cite{Ferrari} where the Laplacian is built from the vector fields of the first layer of the stratification (just as $\mathcal{J}_1$ in the previous example of the Heisenberg group $\mathbb{H}$). \\

In this article we will study this relationship taking into account other types of Laplace operators (also called \textit{sub-Laplacians}) and this will be made in the setting of nilpotent Lie groups that are a generalization of stratified Lie groups. \\

The plan of the article is as follows: first we present the structure of nilpotent Lie groups we are going to work with and then we will state our main theorems. Then in Section \ref{Section_SomeProperties} we introduce the tools that are going to be used and finally in Section \ref{Section_Proofs} we will give the proofs of our results. 
\section{Presentation of the framework and statement of the results}
In this section we will make a short introduction to general nilpotent Lie groups and to some of the basic tools needed in order to present our results. Our main references are the books \cite{Varopoulos} and \cite{Folland2}.\\

Let $\mathbb{G}$ be a connected unimodular Lie group endowed with its Haar measure $dx$.  Denote by $\mathfrak{g}$ the Lie algebra of $\mathbb{G}$ and consider a family (that will be fixed from now on) of left-invariant vector fields on $\mathbb{G}$
\begin{equation}\label{DefinitionFamily}
{\bf X}=\{X_1,...,X_k\}
\end{equation}
satisfying the \textit{Hörmander condition}, which means that the Lie algebra generated by the $X_j$ for $1\leq j\leq k$ is $\mathfrak{g}$. This condition will give several properties of the tools we are going to use here.\\

Our first task is to endow the group $\mathbb{G}$ with a metric structure. Here we have at our disposal the Carnot-Carathéodory metric associated with ${\bf X}$ defined as follows: let $\ell:[0,1]\longrightarrow \mathbb{G}$ be an absolutely continuous path. We say that $\ell$ is admissible if there exists measurable functions $\gamma_1,...,\gamma_k:[0,1]\longrightarrow\mathbb{C}$ such that, for almost every $t\in [0,1]$, we have $\ell'(t)=\sum_{j=1}^k\gamma_j(t)X_j(\ell(t))$. If $\ell$ is admissible, define the length of $\ell$ by $\|\ell\|=\displaystyle{\int_{0}^1}(\sum_{j=1}^k|\gamma_j(t)|^2)^{1/2}dt$. Then, for all $x,y\in \mathbb{G}$, the distance between $x$ and $y$ is the infimum of the lengths of all admissible curves joining $x$ to $y$. We will denote $\|x\|$ the distance between the origin $e$ and $x$ and $\|y^{-1}\cdot x\|$ the distance between $x$ and $y$.\\

Now, we study the behaviour of the Haar measure of the open balls. For $r>0$ and $x\in \mathbb{G}$, denote by $B(x,r)$ the open ball with respect to the Carnot-Carathéodory metric centered in $x$ and of radius $r$, and by 
$$V(r)=\int_{B(x,r)}dx$$
the Haar measure of any ball of radius $r$. When $0<r<1$, there exists $d\in \mathbb{N}^\ast$, $c_l$ and $C_l>0$ such that, for all $0<r<1$ we have
\begin{equation*}
c_l r^d\leq V(r)\leq C_l r^d.
\end{equation*}
The integer $d$ is the \textit{local} dimension of $(\mathbb{G}, {\bf X})$. When $r\geq 1$, two situations may occur, independently of the choice of the family ${\bf X}$: either $\mathbb{G}$ has polynomial volume growth and there exist $D\in \mathbb{N}^\ast$,  $c_\infty$ and $C_\infty>0$ such that, for all $r\geq 1$ we have
\begin{equation*}
c_\infty r^D\leq V(r)\leq C_\infty r^D,
\end{equation*}
or $\mathbb{G}$ has exponential volume growth, which means that there exist  $c_e, C_e, \alpha, \beta>0$ such that, for all $r\geq 1$ we have
$$c_e e^{\alpha r}\leq V(r)\leq C_e  e^{\beta r}.$$
When $\mathbb{G}$ has polynomial volume growth, the integer $D$ is called the dimension at infinity of $\mathbb{G}$.\\

Recall that nilpotent groups have polynomial volume growth and that a \textit{strict} subclass of the nilpotent groups consists of stratified Lie groups where $d=D$. For example, in the case of the Heisenberg group we have $d=D=4$. \\

Once we have fixed the family $\textbf{X}$, we define the gradient on $\mathbb{G}$ by $\nabla = (X_1,...,X_k)$ and we consider a Laplacian $\mathcal{J}$ on $\mathbb{G}$ defined by 
\begin{equation}\label{Definition_LaplacienComplet}
\mathcal{J}=-\sum_{j=1}^k X_j^2,
\end{equation}
which is a positive self-adjoint, hypo-elliptic operator since ${\bf X}$ satisfies the Hörmander's condition, see \cite{Varopoulos}.\\

Its associated heat operator on $]0, +\infty[\times\mathbb{G}$ is given by $\partial_{t}+\mathcal{J}$ and we will denote by $(H_t)_{t> 0}$ the semi-group obtained from the Laplacian $\mathcal{J}$.\\ 

Now, using the spectral theory associated to this Laplacian $\mathcal{J}$, we can define for $0<s<1$ the operator  $\mathcal{J}^s$ (see the details in Section \ref{Section_SomeProperties}).\\

From now on we will work with a unimodular nilpotent Lie group $\mathbb{G}$ with a Haar measure $dx$ and with a family of left-invariant vector fields $\textbf{X}$ that satisfy the Hörmander condition. Here $d$ and $D$ are the local dimension and the dimension at infinity respectively and we will fix the Laplace operator $\mathcal{J}$ given by the formula (\ref{Definition_LaplacienComplet}).\\ 

For more details concerning nilpotent Lie groups see the books \cite{Varopoulos}, \cite{Folland2}, \cite{Stein2} and the articles \cite{Folland0}, \cite{Saka}, \cite{Chamorro} and the references there in.

\begin{Remarque}
The operators $\mathcal{J}_1$ and $\mathcal{J}_3$ given in (\ref{Several_Laplacians}) are in the scope of this framework but the Laplacian $\mathcal{J}_2$ not, as is not bluit as in (\ref{Definition_LaplacienComplet}); however this is not a problem since a small modification of our arguments should lead to similar results in the case of the Heisenberg group $\mathbb{H}$. 
\end{Remarque}

Now that we have fixed our framework we can state our main theorem.
\begin{Theoreme}\label{Theorem_Principal}
Let $\mathbb{G}$ be a nilpotent Lie group and let $\mathcal{J}$ be a Laplace operator. If $\varphi \in \mathcal{S}(\mathbb{G})$ is a smooth function and if we consider the following problem 
\begin{equation}\label{Equation0}
\begin{cases}
\partial^2_t u(t,x)+\frac{1-2s}{t} \partial_t u(t,x)-\mathcal{J}u(t,x)=0,\qquad \mbox{for } x \in  \mathbb{G} \mbox{ and } t>0 \mbox{ with } 0<s<1,\\[5mm]
u(0,x)=\varphi(x).
\end{cases}
\end{equation}
Then the function $u:]0, +\infty[\times \mathbb{G} \longrightarrow \mathbb{R}$ given by the formula
\begin{equation}\label{Def_Solution_Extension_Problem}
u(t,x)=\frac{1}{\Gamma(s)}\int_{0}^{+\infty} H_\tau \mathcal{J}^s \varphi(x)e^{\frac{-t^2}{4\tau}}\frac{d\tau}{\tau^{1-s}}
\end{equation}
is a $L^2$-solution of the extension problem (\ref{Equation0}).\\ 

Moreover, we have the following relationship between the solution of this extension problem and the fractional powers of the Laplacian:
\begin{equation*}
\mathcal{J}^s\varphi(x)=-C(s) \; \underset{t\to 0^+}{\lim}\; t^{1-2s}\partial_t u(t,x).
\end{equation*}
\end{Theoreme}
\section{Some tools and related properties}\label{Section_SomeProperties}
Our proof for the Theorem \ref{Theorem_Principal} relies mainly in the use of two objects and their properties: the heat kernel $h_t$ associated to the Laplacian $\mathcal{J}$ and the spectral decomposition of the Laplacian.\\

We need to fix some terminology. The Lebesgue spaces $L^p(\mathbb{G})$ with $1\leq p\leq +\infty$ are defined for functions $f:\mathbb{G}\longrightarrow \mathbb{R}$ in a classical way by the condition
$$\|f\|_{L^p}=\left(\int_{\mathbb{G}}|f(x)|^p dx\right)^{\frac{1}{p}} \qquad (1\leq p<+\infty)$$
with the usual modification if  $p=+\infty$: $\|f\|_{L^\infty}=\underset{x\in \mathbb{G}}{\sup ess} |f(x)|$. For $f$ and $g$ two measurable functions the convolution $f\ast g$ is given by the expression
$$f\ast g(x)=\int_{\mathbb{G}}f(y)g(y^{-1}\cdot x)dy=\int_{\mathbb{G}}f(x\cdot y^{-1})g(y)dy.$$
Recall that we have at our disposal Young's inequalities (see a proof in \cite{Folland2})
\begin{Lemme}\label{LemYoungG}
Let $1\leq p,q,r\leq +\infty$ such that $1+\frac{1}{r}=\frac{1}{p}+\frac{1}{q}$. If $f\in L^{p}(\mathbb{G})$ and $g\in L^{q}(\mathbb{G})$, then $f\ast g \in L^{r}(\mathbb{G})$ and we have
$$ \|f\ast g\|_{L^r}\leq \|f\|_{L^p}\|g\|_{L^q} $$
\end{Lemme}
Let now $I=(j_{1},...,j_{\beta})\in \{1,...,k\}^{\beta}\; (\beta \in \mathbb{N})$ be a multi-index, we set $|I|=\beta$ and define $X^I$  by the formula $X^{I}=X_{j_{1}}\cdots X_{j_{\beta}}$ with the convention $X^{I}=Id$ if $\beta=0$. The interaction of operators $X^I$  with convolutions is clarified by the following identity:
$$ X^I(f \ast g) = f \ast (X^I g).$$
and we will say that $\varphi\in \mathcal{C}^{\infty}(\mathbb{G})$ belongs to the Schwartz class $\mathcal{S}(\mathbb{G})$ if 
$$N_{\alpha,I}(\varphi)=\underset{x\in \mathbb{G}}{ \sup}(1+\|x\|)^{\alpha}|X^{I}\varphi(x)|<+\infty. \qquad(\alpha\in \mathbb{N}, I\in \underset{\beta\in \mathbb{N}}{\cup} \{1,...,k\}^{\beta}).\\[5mm]$$

One of our fundamental tools is given by the heat semi-group associated to the Laplacian $\mathcal{J}$. We recall in the next theorem some well-known properties of the semi-group $H_t$. 
\begin{Theoreme}\label{fol} There exists a unique family of continuous linear operators $(H_{t})_{t>0}$ defined on $L^{1}+L^{\infty}(\mathbb{G})$ with the semi-group property $H_{t+s}=H_{t}H_{s}$ for all $t, s>0$ and $H_{0}=Id$, such that: 
\begin{enumerate} 
\item[1)] the Laplacian $\mathcal{J}$ is the infinitesimal generator of the semi-group  $H_{t}=e^{-t\mathcal{J}}$;  
\item[2)] $H_{t}$ is a contraction operator on $L^{p}(\mathbb{G})$ for $1\leq p\leq +\infty$ and for $t>0$; 
\item[3)] the semi-group $H_t$ admits a convolution kernel $H_{t}f=f\ast h_{t}$ where $h_{t}$ is the heat kernel.
\item[4)] If $f\in L^{p}(\mathbb{G})$, $1\leq p\leq +\infty$, then the function $u(x, t)=H_{t}f(x)$ is a solution of the heat equation $\partial_tu(t,x)+\mathcal{J}u(t,x)=0$.
\end{enumerate} 
\end{Theoreme} 
We obtain in particular that $H_t$ is a symmetric diffusion semi-group as considered by Stein in \cite{Stein0} with infinitesimal generator $\mathcal{J}$. For a proof and for further details see \cite{Folland2}, \cite{Stein0}, \cite{Varopoulos} and the references given there. \\

The following estimates will be crucial for our purposes. 
\begin{Theoreme}\label{nilpo2}
Let $\mathbb{G}$ be a nilpotent Lie group, $\mathcal{J}$ be the Laplace operator (\ref{Definition_LaplacienComplet}) built from the family  $\textbf{X}$ of vector fields given in (\ref{DefinitionFamily}) and let $h_{t}$ be the associated heat kernel. Then for all $m\in \mathbb{N}$, for all multi-index $I$, there exists a constant $C>0$ such that, for $\varepsilon >0$:
\begin{equation*}
\left|\partial_t^m X^{I} h_{t}(x)  \right|\leq C\, t^{-m-|I|/2}[V(\sqrt{t})]^{-1}exp(-\|x\|^{2}/4(1+\varepsilon)t).
\end{equation*}
for all $x\in \mathbb{G}$ and $t>0$. 
\end{Theoreme}
See a proof of this fact in \cite{Varopoulos}, Theorem \emph{IV.4.2}. A useful inequality is the following one
\begin{Corollaire}\label{nilpo4}
Let $t>0$, $1\leq p\leq +\infty$, $\alpha \in \mathbb{N}$ and $I$ a multi-index. There exists $C>0$ such that:
\begin{equation*}
\left\|(1+\|\cdot\|)^{\alpha}X^{I}h_{t}(\cdot)\right\|_{L^p}\leq C\, t^{-|I|/2}(1+\sqrt{t})^{\alpha}[V(\sqrt{t})]^{-1/p'},\qquad \left(\frac{1}{p}+\frac{1}{p'}=1\right).\\[5mm]
\end{equation*}
\end{Corollaire}

We introduce now our second fundamental tool. Since the Laplacian given by (\ref{Definition_LaplacienComplet}) is a positive self-adjoint operator, it admits a spectral decomposition of the following form
$$\mathcal{J}=\displaystyle{\int_{0}^{+\infty} \lambda dE(\lambda)}$$
This spectral decomposition allows us to define the fractional powers of the Laplacian by the expression
$$\mathcal{J}^s=\displaystyle{\int_{0}^{+\infty} \lambda^s dE(\lambda)}$$
with $0<s$. This formula is very useful to deduce some properties for the fractional Laplacian like $\mathcal{J}^{s_1}[\mathcal{J}^{s_2}f]=\mathcal{J}^{s_1+s_2}f$ whenever these quantities are well defined. But it will also help us to build a family of operators $m(\mathcal{J})$ associated to a Borel function $m$. A classical example of this situation is given by the heat semi-group
$$H_t=e^{-t\mathcal{J}}=\int_{0}^{+\infty} e^{-t\lambda}dE(\lambda)  \qquad \mbox{with } m(\lambda)= e^{-\lambda},$$
and from these formulas we can see that we have the identity $\mathcal{J}^s H_t=H_t \mathcal{J}^s$. It is however possible to go one step further, indeed, following \cite{Hulanicki} and \cite{Furioli2} we have the next result.
\begin{Proposition}\label{Proposition_Crucial}
Let $k\in\mathbb{N}$ and $m$ be a function of class $\mathcal{C}^{k}(\mathbb{R}^{+})$, we write 
\begin{equation*} 
\|m \|_ {(k)}=\underset{\underset{\lambda>0}{0\leq r\leq k}}{\sup } (1+\lambda)^{k }|m^{(r)}(\lambda)|.
\end{equation*}
We define the operator $m(t\mathcal{J})$ for $t>0$ by the expression $\displaystyle{m(t\mathcal{J})=\int_{0}^{+\infty} m(t\lambda) dE(\lambda)}$. Then this operator admits a convolution kernel $M_{t}$ and moreover, for $\alpha\in \mathbb{R}$ and $I$ a multi-index, there exists $C>0$ and $k\in \mathbb{N}$ such that:
\begin{equation*}
\left\|(1+\|\cdot \|)^{\alpha}X^{I}M_{t}(\cdot)\right\|_{L^p}\leq C\, t^{-|I|/2}(1+\sqrt{t})^{\alpha}[V(\sqrt{t})]^{-1/p'}\|m\|_{(k)}\qquad \left(\frac{1}{p}+\frac{1}{p'}=1\right).
\end{equation*}
\end{Proposition}
Note that there will be a size effect in the previous estimate depending on the values of $t$ since we have the polynomial estimates for $V(\sqrt{t})$: more precisely we have $[V(\sqrt{t})]^{-1/p'}\sim t^{^{-d/2p'}}$ if $0<t<1$ and $[V(\sqrt{t})]^{-1/p'}\sim t^{-D/2p'}$, if $t>1$.\\

Now we state a lemma, which is a direct application of the previous result, that will be useful in the sequel:
\begin{Lemme}\label{Lemma_Crucial}
Let $0<s$, then we have the following inequalities
$$\|\mathcal{J}^s h_t\|_{L^1}\leq C t^{-s}\qquad  \mbox{and }\qquad \|\mathcal{J}^s h_t\|_{L^2}\leq C t^{-s}[V(\sqrt{t})]^{-1/2}$$
\end{Lemme}
\textbf{Proof.}
For the $L^1$ estimate we write
\begin{eqnarray*}
\mathcal{J}^s h_t(x)=\mathcal{J}^s(h_{t/2}\ast h_{t/2})(x)&=&\mathcal{J}^s H_{t/2}(h_{t/2})(x)=\left(\int_{0}^{+\infty}\lambda^s e^{-\frac{t}{2}\lambda}dE(\lambda)\right)(h_{t/2})(x)\\
&=&Ct^{-s}\left(\int_{0}^{+\infty}m(t\lambda)dE(\lambda)\right)(h_{t/2})(x)=Ct^{-s} M_t\ast h_{t/2}(x)
\end{eqnarray*}
where $m(\lambda)=\lambda^s e^{-\lambda}$. Then, taking the $L^1$-norm, using Young's inequality and applying the Proposition \ref{Proposition_Crucial} we obtain
$$\|\mathcal{J}^s h_t\|_{L^1}\leq Ct^{-s} \|M_t\|_{L^1}\| h_{t/2}\|_{L^1}=Ct^{-s}\|m\|_{(0)}=Ct^{-s}.$$
The $L^2$ estimate is very similar, indeed, using the identity $\displaystyle{\mathcal{J}^s h_t(x)=Ct^{-s} M_t\ast h_{t/2}(x)}$ obtained above and taking the $L^2$ norm of it we have
$$\|\mathcal{J}^s h_t\|_{L^2}\leq Ct^{-s}\|M_t\|_{L^2}\|h_{t/2}\|_{L^1}=C t^{-s}\|m\|_{(0)}[V(\sqrt{t})]^{-1/2}=C t^{-s}[V(\sqrt{t})]^{-1/2}.$$
\hfill $\blacksquare$

To continue, we will need to introduce one more tool. Following \cite{Stinga1} we will note
\begin{equation}\label{Def_E_fg}
\left\langle \mathcal{J}f,g\right\rangle_{L^{2}}=\int_{0}^{+\infty} \lambda d E_{f,g}(\lambda)
\end{equation}

for all $f \in Dom(\mathcal{J})$ and $g\in L^2(\mathbb{G})$, where $dE_{f,g}$ is regular Borel complex measure of bounded variation concentrated on the spectrum of the operator $\mathcal{J}$.
\section{Proof of the main Theorem}\label{Section_Proofs}
We will follow here some of the ideas of \cite{Stinga1}. Our first task is to proof that each term of the expressions (\ref{Equation0})-(\ref{Def_Solution_Extension_Problem}) is well defined in a $L^2$-sense. 
\begin{Lemme}
Let $\varphi\in \mathcal{S}(\mathbb{G})$ and $t>0$, then the function $u(t,x)$ given by (\ref{Def_Solution_Extension_Problem}) belongs to $L^2(\mathbb{G})$.
\end{Lemme}
\textbf{Proof.}
By duality we consider $g \in L^2(\mathbb{G})$ and we study the quantity 
\begin{eqnarray*}
\langle u(t, \cdot); g\rangle_{L^2}&=&\int_{\mathbb{G}} \left(\frac{1}{\Gamma(s)}\int_{0}^{+\infty} H_\tau \mathcal{J}^s \varphi(x)e^{-\frac{t}{4\tau}}\frac{d\tau}{\tau^{1-s}}\right)g(x)dx\\
&=&\int_{\mathbb{G}} \left(\frac{1}{\Gamma(s)}\int_{0}^{1} H_\tau \mathcal{J}^s \varphi(x)e^{-\frac{t}{4\tau}}\frac{d\tau}{\tau^{1-s}}\right)g(x)dx\\
& &+\int_{\mathbb{G}} \left(\frac{1}{\Gamma(s)}\int_{1}^{+\infty} H_\tau \mathcal{J}^s \varphi(x)e^{-\frac{t}{4\tau}}\frac{d\tau}{\tau^{1-s}}\right)g(x)dx
\end{eqnarray*}
and applying the properties of the heat kernel and of the fractional powers of the Laplacian we write
\begin{eqnarray*}
\langle u(t, \cdot); g\rangle_{L^2}&=&\int_{\mathbb{G}} \left(\frac{1}{\Gamma(s)}\int_{0}^{1} \mathcal{J}^{\frac{s}{2}}\varphi  \ast \mathcal{J}^{\frac{s}{2}} h_\tau (x)e^{-\frac{t}{4\tau}}\frac{d\tau}{\tau^{1-s}}\right)g(x)dx\\
&+&\int_{\mathbb{G}} \left(\frac{1}{\Gamma(s)}\int_{1}^{+\infty} \varphi  \ast \mathcal{J}^{s} h_\tau (x)e^{-\frac{t}{4\tau}}\frac{d\tau}{\tau^{1-s}}\right)g(x)dx\\
\end{eqnarray*}
Thus, by the Cauchy-Schwarz inequality and by Young's inequality we obtain
\begin{eqnarray*}
\left|\langle u(t, \cdot); g\rangle_{L^2}\right|& \leq & \frac{1}{\Gamma(s)}\|g\|_{L^2}\left( \int_{0}^{1}\|\mathcal{J}^{\frac{s}{2}}\varphi\|_{L^2}\|\mathcal{J}^{\frac{s}{2}}h_\tau\|_{L^1}e^{-\frac{t}{4\tau}}\frac{d\tau}{\tau^{1-s}}+\int_{1}^{+\infty}\|\varphi\|_{L^1}\|\mathcal{J}^{s} h_\tau\|_{L^2}e^{-\frac{t}{4\tau}}\frac{d\tau}{\tau^{1-s}}\right)
\end{eqnarray*}
Now we use the estimates available in the Lemma \ref{Lemma_Crucial} for the heat kernel:
\begin{eqnarray*}
\left|\langle u(t, \cdot); g\rangle_{L^2}\right|&\leq & \frac{1}{\Gamma(s)}\|g\|_{L^2}\left( C\int_{0}^{1}\|\mathcal{J}^{\frac{s}{2}}\varphi\|_{L^2}  \tau^{-\frac{s}{2}} \; e^{-\frac{t}{4\tau}}\frac{d\tau}{\tau^{1-s}}+C'\int_{1}^{+\infty}\|\varphi\|_{L^1} \tau^{-s} [V(\sqrt{\tau})]^{-\frac{1}{2}} e^{-\frac{t}{4\tau}}\frac{d\tau}{\tau^{1-s}}\right)\\
&\leq & \frac{C''}{\Gamma(s)}\|g\|_{L^2}\left(\|\mathcal{J}^{\frac{s}{2}}\varphi\|_{L^2} \int_{0}^{1}  \tau^{\frac{s}{2}-1}d\tau+\|\varphi\|_{L^1} \int_{1}^{+\infty}\tau^{-\frac{D}{4}-1} d\tau\right)\\
&\leq & \frac{C}{\Gamma(s)}\|g\|_{L^2} \left(\|\mathcal{J}^{\frac{s}{2}}\varphi\|_{L^2}+\|\varphi\|_{L^1}\right).
\end{eqnarray*}
and we obtain that the function $u(t,x)$ is in $L^2(\mathbb{G})$ for all $t>0$. 
\hfill $\blacksquare$

\begin{Lemme}
$u(t,x)$ belongs to the domain of the operator $\mathcal{J}$. 
\end{Lemme}
\textbf{Proof.} We need to prove that for all $g\in L^2(\mathbb{G})$ we have
$$ \underset{\rho \to 0^+}{\lim} \left\langle \frac{H_\rho u(t, \cdot)-u(t, \cdot) }{\rho}; g \right\rangle_{L^2}<+\infty$$
so we study the following quantity
\begin{eqnarray*}
\left\langle \frac{H_\rho u(t, \cdot)-u(t, \cdot) }{\rho}; g \right\rangle_{L^2}&=&\int_{\mathbb{G}}\frac{H_\rho u(t, x)-u(t, x) }{\rho} g(x)dx
\end{eqnarray*}
\begin{eqnarray*}
&=&\frac{1}{\Gamma(s)}\int_{\mathbb{G}}\frac{1}{\rho}\left(H_\rho\left(\int_{0}^{+\infty} H_\tau \mathcal{J}^s \varphi(x)e^{-\frac{t^2}{4\tau}}\frac{d\tau}{\tau^{1-s}}\right)-\int_{0}^{+\infty} H_\tau \mathcal{J}^s \varphi(x)e^{-\frac{t^2}{4\tau}}\frac{d\tau}{\tau^{1-s}} \right) g(x)dx\\
&=& \frac{1}{\Gamma(s)}\int_{\mathbb{G}} \int_{0}^{+\infty} \left(\frac{H_{\rho+\tau}\mathcal{J}^s\varphi(x)-H_\tau \mathcal{J}^s\varphi(x)}{\rho}\right) e^{-\frac{t^2}{4\tau}}\frac{d\tau}{\tau^{1-s}} \; g(x)dx
\end{eqnarray*}
Taking $\rho \to 0^+$ we have
\begin{eqnarray*}
\underset{\rho \to 0^+}{\lim} \left\langle \frac{H_\rho u(t, \cdot)-u(t, \cdot) }{\rho}; g \right\rangle_{L^2}&=&\frac{1}{\Gamma(s)}\int_{\mathbb{G}} \int_{0}^{+\infty} \mathcal{J}H_\tau \mathcal{J}^s\varphi(x) e^{-\frac{t^2}{4\tau}}\frac{d\tau}{\tau^{1-s}} \; g(x)dx\\
&\leq & \frac{1}{\Gamma(s)}\int_{0}^{+\infty} \|\mathcal{J}^{s+1}H_\tau \varphi\|_{L^2} \|g\|_{L^2} e^{-\frac{t^2}{4\tau}}\frac{d\tau}{\tau^{1-s}}\\
&\leq & \frac{\|g\|_{L^2}}{\Gamma(s)}\left(\int_{0}^1 \|\mathcal{J}^{s+1}H_\tau \varphi\|_{L^2} e^{-\frac{t^2}{4\tau}}\frac{d\tau}{\tau^{1-s}}+\int_{1}^{+\infty} \|\mathcal{J}^{s+1}H_\tau \varphi\|_{L^2} e^{-\frac{t^2}{4\tau}}\frac{d\tau}{\tau^{1-s}}\right)
\end{eqnarray*}
Now, using the properties of the heat kernel and the special relationships of this kernel with the fractional powers of the Laplacian we can write
\begin{eqnarray*}
&\leq & \frac{\|g\|_{L^2}}{\Gamma(s)}\left(\int_{0}^1 \|\mathcal{J}^{\frac{s}{2}+1}\varphi \ast \mathcal{J}^{\frac{s}{2}}h_\tau\|_{L^2} e^{-\frac{t^2}{4\tau}}\frac{d\tau}{\tau^{1-s}}+\int_{1}^{+\infty} \| \varphi\ast \mathcal{J}^{s+1}h_\tau\|_{L^2} e^{-\frac{t^2}{4\tau}}\frac{d\tau}{\tau^{1-s}}\right)\\
&\leq & \frac{\|g\|_{L^2}}{\Gamma(s)}\left(\int_{0}^1 \|\mathcal{J}^{\frac{s}{2}+1}\varphi\|_{L^2}\|\mathcal{J}^{\frac{s}{2}}h_\tau\|_{L^1} e^{-\frac{t^2}{4\tau}}\frac{d\tau}{\tau^{1-s}}+\int_{1}^{+\infty} \|\varphi\|_{L^2} \|\mathcal{J}^{s+1}h_\tau\|_{L^1} e^{-\frac{t^2}{4\tau}}\frac{d\tau}{\tau^{1-s}}\right)
\end{eqnarray*}
Now, using inequalities stated in Lemma \ref{Lemma_Crucial} for the heat kernel we obtain
\begin{eqnarray*}
\underset{\rho \to 0^+}{\lim} \left\langle \frac{H_\rho u(t, \cdot)-u(t, \cdot) }{\rho}; g \right\rangle_{L^2}&\leq & \frac{\|g\|_{L^2}}{\Gamma(s)}\left(C_1\int_{0}^1 \|\mathcal{J}^{\frac{s}{2}+1}\varphi\|_{L^2}\tau^{\frac{s}{2}-1} d\tau+C_2\int_{1}^{+\infty} \|\varphi\|_{L^2} \tau^{-2}d\tau\right)\\
&\leq & \frac{C\|g\|_{L^2}}{\Gamma(s)}\left(\|\mathcal{J}^{\frac{s}{2}+1}\varphi\|_{L^2}+\|\varphi\|_{L^2} \right).
\end{eqnarray*}
\hfill $\blacksquare$\\

\begin{Lemme}
We have $\displaystyle{\underset{t \to 0^+}{\lim} u(t,x)=\varphi(x)}$.
\end{Lemme}
\textbf{Proof.} For a function $g\in L^2(\mathbb{G})$ we have
\begin{eqnarray*}
\langle u(t,\cdot); g \rangle_{L^2}&=&\int_{\mathbb{G}} \left(\frac{1}{\Gamma(s)}\int_{0}^{+\infty} H_\tau \mathcal{J}^s \varphi(x)e^{-\frac{t^2}{4\tau}}\frac{d\tau}{\tau^{1-s}} \right) g(x)dx\\
&=&\frac{1}{\Gamma(s)}\int_{0}^{+\infty}\left(\int_{\mathbb{G}} H_\tau \mathcal{J}^s \varphi(x)g(x) dx\right) e^{-\frac{t^2}{4\tau}}\frac{d\tau}{\tau^{1-s}}\\
&=& \frac{1}{\Gamma(s)}\int_{0}^{+\infty}\left(\int_{\mathbb{G}} \left(\int_{0}^{+\infty}e^{-\tau \lambda}\lambda^s dE_\lambda\right) \varphi(x)g(x) dx\right) e^{-\frac{t^2}{4\tau}}\frac{d\tau}{\tau^{1-s}}
\end{eqnarray*}
by definition of the measure $d E_{\varphi, g}$ given in (\ref{Def_E_fg}) and by Fubini's theorem we can write

\begin{eqnarray*}
\langle u(t,\cdot); g \rangle_{L^2}&=&\frac{1}{\Gamma(s)}\int_{0}^{+\infty}\int_{0}^{+\infty}e^{-\tau \lambda}\lambda^s dE_{\varphi, g}(\lambda)  e^{-\frac{t^2}{4\tau}}\frac{d\tau}{\tau^{1-s}}\\
&=&\frac{1}{\Gamma(s)}\int_{0}^{+\infty}\int_{0}^{+\infty}e^{-\tau \lambda}(\tau\lambda)^s  e^{-\frac{t^2}{4\tau}} \frac{d\tau}{\tau} dE_{\varphi, g}(\lambda)
\end{eqnarray*}
now, if we set $\rho= \lambda \tau$ we obtain
\begin{eqnarray*}
\langle u(t,\cdot); g \rangle_{L^2}=\frac{1}{\Gamma(s)}\int_{0}^{+\infty}\int_{0}^{+\infty}e^{-\rho}\rho^s  e^{-\frac{\lambda t^2}{4\rho}} \frac{d\rho}{\rho} dE_{\varphi, g}(\lambda)
\end{eqnarray*}
At this point we take the limit $t\to 0^+$:
\begin{eqnarray*}
\underset{t\to 0^+}{\lim}\langle u(t,\cdot); g \rangle_{L^2}&=&\underset{t\to 0^+}{\lim}\frac{1}{\Gamma(s)}\int_{0}^{+\infty}\int_{0}^{+\infty}e^{-\rho}\rho^s  e^{-\frac{\lambda t^2}{4\rho}} \frac{d\rho}{\rho} dE_{\varphi, g}(\lambda)\\
&=&\frac{1}{\Gamma(s)}\int_{0}^{+\infty}\int_{0}^{+\infty}e^{-\rho}\rho^{s-1} d\rho \;dE_{\varphi, g}(\lambda)\\
&=& \int_{0}^{+\infty} dE_{\varphi, g}(\lambda)=\langle \varphi ; g \rangle_{L^2}
\end{eqnarray*}
\hfill $\blacksquare$\\

\begin{Lemme}
The function $u(t,x)$ is derivable with respect to $t$ and we have the following formula
$$\partial_t u(t,x)=\frac{-1}{\Gamma(s)}\int_{0}^{+\infty} H_\tau \mathcal{J}^s\varphi(x) \frac{t e^{-\frac{t^2}{4\tau}}}{2\tau }   \frac{d\tau}{\tau^{1-s}}$$
\end{Lemme}
\textbf{Proof.} Our starting point is the expression
\begin{eqnarray*}
\underset{h\to 0^+}{\lim} \left\langle \frac{u(t+h,\cdot)-u(t, \cdot)}{h}; g\right \rangle_{L^2}&=&\frac{1}{\Gamma(s)}\int_{\mathbb{G}}\int_{0}^{+\infty} H_\tau \mathcal{J}^s\varphi(x)\;\left( \underset{h\to 0^+}{\lim}\frac{e^{-\frac{(t+h)^2}{4\tau}}-e^{-\frac{t^2}{4\tau}} }{h}\right) \frac{d\tau}{\tau^{1-s}} g(x)dx\\
&=& \frac{1}{\Gamma(s)}\int_{\mathbb{G}}\int_{0}^{+\infty} H_\tau \mathcal{J}^s\varphi(x)\;\partial_t(e^{-\frac{t^2}{4\tau}}) \frac{d\tau}{\tau^{1-s}} g(x)dx\\
&=& \int_{\mathbb{G}}\left(\frac{-1}{\Gamma(s)}\int_{0}^{+\infty} H_\tau \mathcal{J}^s\varphi(x) \frac{t e^{-\frac{t^2}{4\tau}}}{2\tau }   \frac{d\tau}{\tau^{1-s}}\right) g(x)dx.
\end{eqnarray*}
\hfill $\blacksquare$\\
\begin{Corollaire}
The function $u(t,x)$ is twice differenciable with respect to $t$ and we have the formula
$$\partial_t^2 u(t,x)=\frac{1}{\Gamma(s)}\int_{0}^{+\infty}H_\tau \mathcal{J}^s\varphi(x)\left(\frac{t^2}{4\tau^2}-\frac{1}{2\tau} \right)e^{-\frac{t^2}{4\tau}} \frac{d\tau}{\tau^{1-s}}$$
\end{Corollaire}
Now that we have verified that all the terms of the extension problem (\ref{Equation0}) are well defined in the $L^2$-sense, we can state the following proposition
\begin{Proposition}
The function $u(t,x)$ defined by (\ref{Def_Solution_Extension_Problem}) and associated to an initial data $\varphi\in \mathcal{S}(\mathbb{G})$ satisfies in the $L^2$-sense the equation
$$\partial_t^2 u(t,x)+\frac{1-2s}{t}\partial_t u(t,x)-\mathcal{J}u(t,x)=0$$
\end{Proposition}
\textbf{Proof.} Let $g\in L^2(\mathbb{G})$, we have
$$\left\langle \partial_t^2 u(t,\cdot)+\frac{1-2s}{t}\partial_t u(t,\cdot)-\mathcal{J}u(t,\cdot); g \right\rangle_{L^2}= \int_{\mathbb{G}} \frac{1}{\Gamma(s)}\int_{0}^{+\infty}H_\tau \mathcal{J}^s\varphi(x)\left(\frac{t^2}{4\tau^2}-\frac{1}{2\tau} \right)e^{-\frac{t^2}{4\tau}} \frac{d\tau}{\tau^{1-s}} g(x)dx$$
$$+ \int_{\mathbb{G}} \frac{1-2s}{t} \left(\frac{-1}{\Gamma(s)}\int_{0}^{+\infty} H_\tau \mathcal{J}^s\varphi(x) \frac{t e^{-\frac{t^2}{4\tau}}}{2\tau }   \frac{d\tau}{\tau^{1-s}}\right) g(x)dx - \int_{\mathbb{G}}\mathcal{J}u(t,x)g(x)dx$$
\begin{eqnarray*}
&=&\int_{\mathbb{G}}\left(\frac{1}{\Gamma(s)}\int_{0}^{+\infty}  H_\tau \mathcal{J}^s\varphi(x)\left(\frac{t^2}{4\tau^2}+\frac{s-1}{\tau} \right)e^{-\frac{t^2}{4\tau}} \frac{d\tau}{\tau^{1-s}}\right)g(x)dx-\int_{\mathbb{G}}\mathcal{J}u(t,x)g(x)dx
\end{eqnarray*}
At this point, we can perform an integration by parts in the first integral with respect to the variable $\tau$ to obtain 
\begin{eqnarray*}
\left\langle \partial_t^2 u(t,\cdot)+\frac{1-2s}{t}\partial_t u(t,\cdot)+\mathcal{J}u(t,\cdot); g \right\rangle_{L^2}&=&
\frac{-1}{\Gamma(s)}\int_{\mathbb{G}} \int_{0}^{+\infty} \partial_\tau \left(H_\tau \mathcal{J}^s\varphi(x)\right)e^{-\frac{t^2}{4\tau}} \frac{d\tau}{\tau^{1-s}}\\
& &-\int_{\mathbb{G}}\mathcal{J}u(t,x)g(x)dx
\end{eqnarray*}
But, since $\partial_\tau H_\tau f=-\mathcal{J}H_\tau f$ we obtain
\begin{eqnarray*}
\left\langle \partial_t^2 u(t,\cdot)+\frac{1-2s}{t}\partial_t u(t,\cdot)+\mathcal{J}u(t,\cdot); g \right\rangle_{L^2}&=&
\frac{1}{\Gamma(s)}\int_{\mathbb{G}} \int_{0}^{+\infty}  \mathcal{J}\left( H_\tau \mathcal{J}^s\varphi(x)\right) e^{-\frac{t^2}{4\tau}} \frac{d\tau}{\tau^{1-s}}\\
& &-\int_{\mathbb{G}}\mathcal{J}u(t,x)g(x)dx\\
&=& \int_{\mathbb{G}}\mathcal{J}u(t,x)g(x)dx-\int_{\mathbb{G}}\mathcal{J}u(t,x)g(x)dx=0.
\end{eqnarray*}
This proof the first part of the Theorem \ref{Theorem_Principal}. \hfill $\blacksquare$\\

It remains to proof the last part of our main theorem and this is done with the next proposition. 
\begin{Proposition} Let $\varphi \in \mathcal{S}(\mathbb{G})$ be an initial data of the extension problem (\ref{Equation0}), let $u(t,x)$ be the function defined by the formula (\ref{Def_Solution_Extension_Problem}), then in the $L^2$-sense we have the limit
$$\underset{t\to 0^+}{ \lim} t^{1-2s}\partial_tu(t,x)=-C(s)\mathcal{J}^s\varphi(x).$$
\end{Proposition} 
\textbf{Proof.} For $g\in L^2(\mathbb{G})$ we have
\begin{eqnarray*}
\langle  t^{1-2s}\partial_tu(t,\cdot); g \rangle_{L^2}=\int_{\mathbb{G}} t^{1-2s}\partial_tu(t,x)g(x)dx=\frac{-1}{\Gamma(s)} \int_{\mathbb{G}} \int_{0}^{+\infty}  H_\tau\mathcal{J}^s\varphi(x)  \frac{t^{2-2s}}{2\tau^{2-s}}e^{-\frac{t^2}{4\tau}}d\tau g(x)dx
\end{eqnarray*}
Making the change of variables $u=\frac{\tau}{t^2}$ we obtain
\begin{eqnarray*}
\langle  t^{1-2s}\partial_tu(t,\cdot); g \rangle_{L^2}= \frac{-1}{\Gamma(s)} \int_{\mathbb{G}} \int_{0}^{+\infty}  H_{t^2u}\mathcal{J}^s\varphi(x)  \frac{ e^{-\frac{1}{4 u}}}{2u^{2-s}}d u\; g(x)dx.
\end{eqnarray*}
Now, taking $t\to 0^+$ we have
\begin{eqnarray*}
\underset{t\to 0^+}{ \lim} \langle  t^{1-2s}\partial_tu(t,\cdot); g \rangle_{L^2}&=&\underset{t\to 0^+}{ \lim} \frac{-1}{\Gamma(s)} \int_{\mathbb{G}} \int_{0}^{+\infty}  H_{t^2u}\mathcal{J}^s\varphi(x)  \frac{ e^{-\frac{1}{4 u}}}{2u^{2-s}}d u\; g(x)dx\\
&=& \frac{-1}{\Gamma(s)} \int_{\mathbb{G}} \int_{0}^{+\infty} \mathcal{J}^s\varphi(x)  \frac{ e^{-\frac{1}{4 u}}}{2u^{2-s}}d u\; g(x)dx\\
&=&\frac{-1}{\Gamma(s)} \int_{\mathbb{G}}\mathcal{J}^s\varphi(x)   \left(\int_{0}^{+\infty} \frac{ e^{-\frac{1}{4 u}}}{2u^{2-s}}d u\right)\; g(x)dx\\
&=& -C(s)\int_{\mathbb{G}}\mathcal{J}^s\varphi(x) g(x)dx.
\end{eqnarray*}
\hfill $\blacksquare$

\begin{Remarque}
Nilpotent Lie groups have some special properties compared to general polynomial volume growth Lie groups, see the details in \cite{Varopoulos}. However, all the properties used here for the heat kernel and for the spectral decomposition for a Laplacian satisfying the Hörmander condition remains true in the general setting of polynomial volume growth Lie groups (see again the book \cite{Varopoulos} and Section 4 of \cite{Chamorro}). Theorem \ref{Theorem_Principal} is still valid in this general setting as the proof follows the same lines exposed here. 
\end{Remarque}


\quad\\[5mm]

\begin{flushright}
\begin{minipage}[r]{80mm}
Diego \textsc{Chamorro} \& Oscar \textsc{Jarr\'in}\\[3mm]
Laboratoire de Math\'ematiques\\ et Mod\'elisation d'Evry\\ (LaMME) UMR 8071\\
Université d'Evry Val d'Essonne\\[2mm]
23 Boulevard de France\\
91037 Evry Cedex\\[2mm]
diego.chamorro@univ-evry.fr
\end{minipage}
\end{flushright}


\begin{thebibliography}{2}
\bibitem{Banica}
V. \textsc{Banica}, M.d.M. \textsc{Gonz\'alez}, M. \textsc{S\'aez}, \emph{Some constructions for the fractional Laplacian on noncompact manifolds}, to appear in Rev. Mat. Iberoam
\bibitem{Cabre}
X. \textsc{Cabr\'e} \&  Y. \textsc{Sire}, \emph{Non linear equations for fractional Laplacians I: regularity, maximum principles and Hamiltonian estimates}. Annales de l'Institut Henri Poincare (C) Non Linear Analysis Volume 31, Issue 1, January–February (2014), Pages 23–53.
\bibitem{Caffarelli}
L. \textsc{Caffarelli} \& L. \textsc{Silvestre}, \emph{An Extention Problem Related to the Fractional Laplacian}, Communications in Partial Differential Equations. Vol 32: 1245-1260, (2007).
\bibitem{Chamorro}
D. \textsc{Chamorro}, \emph{Some functional inequalities on polynomial volume growth Lie groups}, Canad. J. Math. 64 (2012), 481-496.
\bibitem{Ferrari}
F. \textsc{Ferrari} \& B. \textsc{Franchi}, \emph{Hanarck inequality for fractional sub-Laplacians in Carnot groups}, arxiv.org/pdf/1206.0885
\bibitem{Folland0}
G. \textsc{Folland}, \emph{Subelliptic estimates and function spaces on nilpotent Lie groups}.
Ark. Mat. 13 (1975), p. 161-208.
\bibitem{Folland2}
G. \textsc{Folland} \& E. M. \textsc{Stein}, \emph{Hardy Spaces on homogeneous groups}. Mathematical Notes, 28, Princeton University Press (1982).
\bibitem{Frank}
R. \textsc{Frank}, M.D.M. \textsc{Gonz\'alez}, D. D. \textsc{Monticelli} \& J. \textsc{Tan}, \emph{An extension problem for the CR fractional Laplacian}. (arXiv:1312.3381)
\bibitem{Furioli1}
G. \textsc{Furioli}, C. \textsc{Melzi} \& A. \textsc{Veneruso}, \emph{Strichartz inequalities for the wave equation with the full Laplacien on the Heisenberg group}, Studia Math. 160 (2004), 157-178.
\bibitem{Furioli2}
G. \textsc{Furioli}, C. \textsc{Melzi} \& A. \textsc{Veneruso}, \emph{Littlewood-Paley decomposition and Besov spaces on Lie groups of polynomial growth}, Mathematische Nachrichten, vol. 279, no. 9-10, pp. 1028–1040, (2006).
\bibitem{Gale}
J. E. \textsc{Gal\'e}, P. J. \textsc{Miana} \& P. \textsc{Stinga}, \emph{Extension problem and fractional operators: semigroups and wave equations}, J. Evol. Equ. 13 (2013), 343–368.
\bibitem{Hulanicki}
A. \textsc{Hulanicki}, \emph{A functional calculus for Rockland operators on nilpotent Lie groups}. Studia Mathematica T. LXXVIII, (1984).
\bibitem{Saka}
K.\textsc{Saka}, \emph{Besov Spaces and Sobolev spaces on a nilpotent Lie group}. Tohoku. Math. Journ.  Vol. 31, (1979), p. 383-437.
\bibitem{Stein0}
E. M. \textsc{Stein}, \emph{Topics in Harmonic analysis}. Annals of mathematics studies, 63.
Princeton University Press (1970).
\bibitem{Stein2}
E. M. \textsc{Stein}. \emph{Harmonic Analysis}, Princeton University Press (1993).
\bibitem{Stinga1}
P. \textsc{Stinga}, J. \textsc{Torrea}. \emph{Extension Problem and Harnack's inequality for some fractional operators}, Communications in Partial Differential Equations Volume 35, Issue 11, (2010).
\bibitem{Varopoulos}
N. Th. \textsc{Varopoulos}, L. \textsc{Saloff-Coste} \& T. \textsc{Coulhon}, \emph{Analysis and geometry on groups}.
Cambridge Tracts in Mathematics, Vol. 100. (1992).
\end{thebibliography}
\end{document}